\documentclass[]{llncs}      
\usepackage[dvips]{graphicx} 
\begin{document}             

\title{Using matrices  \\ 
		   in post-processing phase \\
			 of CFD simulations}

\author{\bf Gianluca Argentini}
\institute{
         {\it Advanced Computing Laboratory} \\
         Information and Communication Technology Department \\
         RIELLO Group \\
         37045 Legnago (Verona), Italy \\
         gianluca.argentini@riellogroup.com}
        
\maketitle

\begin{abstract}
  In this work I present a technique of construction and fast 
evaluation of a family of cubic polynomials for analytic smoothing 
and graphical rendering of particles trajectories for flows in a 
generic geometry. The principal result of the work was implementation
and test of a method for interpolating 3D points by regular 
parametric curves and their fast and efficient evaluation for a good 
resolution of rendering. For the purpose I have used a parallel 
environment using a multiprocessor cluster architecture.
The efficiency of the used method is good, mainly reducing the number 
of floating-points computations by caching the numerical values of some
line-parameter's powers, and reducing the necessity of communication 
among processes. This work has been developed for the Research \& 
Development Department of my company for planning advanced customized 
models of industrial burners. 
\end{abstract}

\section{Introduction}
  Industrial and power burners have some particular requirements, as a
customized study of the geometry for combustion head and combustion
chamber for an optimal shape of the flame. Rapid prototyping for an
accurate design of the correct geometry involves a numerical simulation
of the gas or oil flows in the burner's components (see Fig. 1).

\begin{figure}[]
\begin{center}
\includegraphics{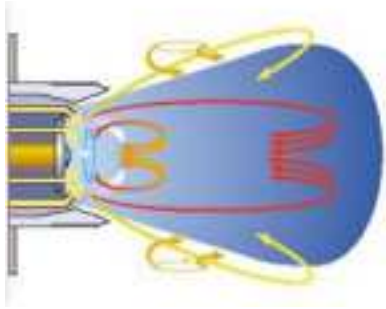}
\end{center}
\caption{Flows in combustion head (left) and combustion chamber.}\label{fig:Head1}
\end{figure}

The necessity of an high graphic resolution require a large amount of 
particles paths for tracing the streamlines of flow. Hence the numerical 
computation is memory and cpu very expensive for the used hardware 
environment. In a tipical simulation the number of paths to compute is
some thousands, and the number of geometrical points to interpolate for
each path is some thousands too. For the treatment of this large amount 
of data a parallel environment can be very useful.

\section{Fitting trajectories with cubic polynomials}
We suppose to have a dataset output from pre-processing and processing 
phases of a simulation, for example from numerical resolution of 
Navier-Stokes equations or from Cellular Automaton models [1]. We would 
a fast and flexible method to obtain from those data an accurate paths
tracking of fluid particles with a smooth 3D visualization of trajectories, 
possibly with continuous slope and curvature. Our experience 
shows that Computational Fluid Dynamics packages have some
limits in this post-processing phase, principally due to a rigid resolution
of the initial mesh and to a small degree of parallelism.

Let {\bf S} the number of 3D points for each trajectory and {\bf M} the
total number of trajectories from simulation dataset. We have tested
that usual interpolation methods are have some disadvantages for
our aims: for example Bezier-like is not realistic in case of twisting
or diverging speed-fields; Chebychev or Least-Squares-like are too
rigid for a customized application; polynomial fitting is simple but
often shows spurious effects as Runge-Gibbs phenomenon [2]. We have
elaborated a {\it spline}-based technique.

We suppose {\bf S} = 4x{\bf N}. For every group of four points, the 
interpolation is obtained by three cubic polynomials imposing four
analytical conditions: passage at {\bf P}$_k$ point, $1\leq{k}\leq{3}$;
passage at {\bf P}$_{k+1}$ point; continuous slope at {\bf P}$_k$ point;
continuous curvature at {\bf P}$_k$ point. For smooth rendering and
for avoiding excessive twisting of trajectories, the cubics {\bf u}$_k$
are added to the Bezier curve {\bf b} associated to the four points:
{\bf v}$_k$ = $\alpha${\bf b} + $\beta${\bf u}$_k$, $0<\alpha,\beta<1$ 
(see Fig. 2).

\begin{figure}[]
\begin{center}
\includegraphics{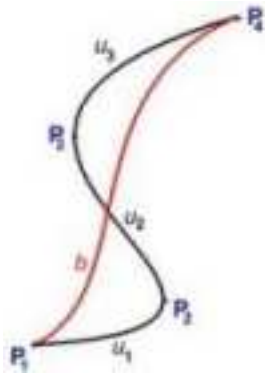}
\end{center}
\caption{Spline-based method with continuous slope and curvature.}\label{fig:Cubics1}
\end{figure}

In our simulations we have chosen $\alpha$ = $\beta$ = 0.5 .
Let {\bf b} = $As^3+Bs^2+Cs+D$, $0\leq{s}\leq{1}$, the Bezier curve
of control points {\bf P}$_1$, ..., {\bf P}$_4$, and let {\bf u}$_k$ 
= $at^3+bt^2+ct+d$, $0\leq{t}\leq{1}$, the spline between two points.
One can see that the coefficients of this spline can be computed by
a matrix-vector product {\bf coeff} = {\bf T}*{\bf p} where {\bf coeff} =
($a$, $b$, $c$, $d$), {\bf p} = ({\bf P}$_{k+1}$, {\bf P}$_k$, $B$, $C$, 1) and
{\bf T} is a 4 x 5 numerical matrix, constant for every groups of points
and for every trajectory. If we define the 4{\bf M}x5{\bf M} 
{\em global matrix}

\[ {\bf G} = \left( \begin{array}{ccccc}
{\bf T} & {\bf 0} & . & . & {\bf 0} \\
{\bf 0} & {\bf T} & . & . & {\bf 0} \\
. & . & . & . & . \\
. & . & . & . & . \\
{\bf 0} & {\bf 0} & . & . & {\bf T} \end{array} \right)\]
where {\bf 0} is a 4x5 zero-matrix, and define the vector {\bf s} = 
({\bf P}$_{k+1}$, {\bf P}$_k$, $B_1$, $C_1$, 1, ..., {\bf P}$_{k+1}$, 
{\bf P}$_k$, $B_{\bf M}$, $C_{\bf M}$, 1), one can compute for every two-points
group the coefficients of cubic splines for all the {\bf M} trajectories with
the matrix-vector product {\bf c} = {\bf G}*{\bf s}. The matrix {\bf G} is
{\em sparse} with density equal at most to 1/{\bf M}; if {\bf M} = 1000, 
the density of 0.001 is a very good value for obtain the benefits of 
sparsity methods, mainly in computational total time and memory allocation [3].

\section{Computing splines}
For computing the coefficients of all the splines involved in the
simulation, the complexity analysis shows a total number of operations of
order {\bf M}*{\bf N}. Using {\bf P} computational processes on a 
multiprocessor environment, a useful method is the distribution of 
{\bf M}/{\bf P} trajectories to every process. In this way every process
receives {\bf M}/{\bf P} rows of the matrix {\bf G} for computing splines
by matrix-vector multiply.
In a first experiment (fall 2003), we have used the Linux cluster at CINECA,
Bologna (Italy), equipped with Pentium III 1.133 GHz processors, and a
software environment constituted by C programs and MPI libraries [4]. The
use of such parallel routines has been useful only for startup of multi-
processes and data distribution. Tests have shown a quasi-linear {\em speedup}, 
in the sense of parallelism, for all the values of {\bf M} and {\bf N}
respect to the number {\bf P} of used processes (see e.g. Fig. 3).

\begin{figure}[]
\begin{center}
\includegraphics[height=5.0cm,width=7.0cm]{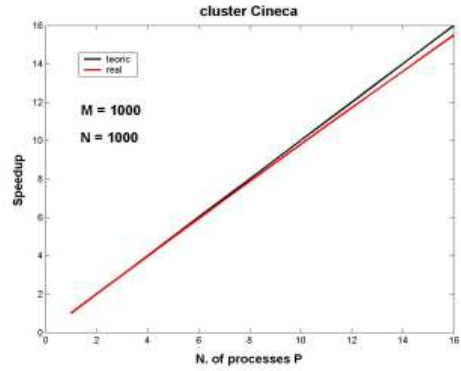}
\end{center}
\caption{Speedup (red line) in computing splines coefficients using C routines and
a Linux cluster.}\label{fig:Cineca1}
\end{figure}

In a second experiment (winter 2003), we have used a multinode Windows 2000 cluster
of our company, equipped with a total of 4 Intel Xeon 3.2 GHz processors and 4 GB Ram, 
and a parallel environment using MATLAB 6.5 scripts on distributed package's 
sessions on nodes. Tests have shown very high performances for splines computation
using the internal algorithms of sparse matrix-vector multiply for the
matrix {\bf G} (see e.g. Fig. 4).

\begin{figure}[]
\begin{center}
\includegraphics[height=7.0cm,width=8.0cm]{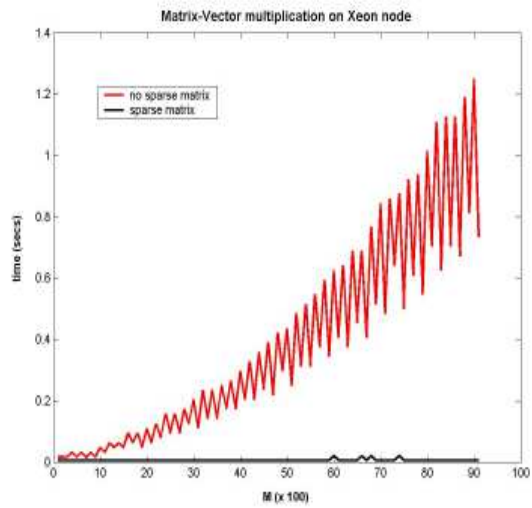}
\end{center}
\caption{Total CPU-time for {\bf G}*{\bf s} matrix-vector multiplication on
a Xeon 3.2 GHz node, using sparse (black line) and full (red line) matrix Matlab algorithm.}\label{fig:Sparse1}
\end{figure}

\section{Valuating splines}
After the computation of splines, we have focused on their valuations on a
suitable set of parameter's values. This set can be chosen large enough to
obtain a fine sampling for an high graphic resolution. Consequently the 
amount of computation can be very huge, so that it is necessary an 
adequate method to valuate all the splines for all the trajectories.

Let {\bf V}+1 the number of ticks for each spline valuation with a uniform
sampling; then the ticks are (0, 1/{\bf V}, 2/{\bf V}, . . ., ({\bf V}-1)/{\bf V}, 1).
The values of splines parameter {\em t} are (0, 1, 2, 3)-th degree powers of this
array. The value of a cubic at {\em t}$_0$ can be view as a dot product:

$at_0^3+bt_0^2+ct_0+d = (a, b, c, d) . (t_0^3, t_0^2, t_0, 1)$

This fact permits to consider the constant 4x({\bf V}+1) matrix

\[ {\bf T} = \left( \begin{array}{cccccc}
0 & (1/{\bf V})^3 & . & . & (({\bf V}-1)/{\bf V})^3 & 1 \\
0 & (1/{\bf V})^2 & . & . & (({\bf V}-1)/{\bf V})^2 & 1 \\
0 & (1/{\bf V})^1 & . & . & (({\bf V}-1)/{\bf V})^1 & 1 \\
1 & 1 & . & . & 1 & 1 \end{array} \right)\]

We consider the {\bf M}x4 matrix 

\[ {\bf C} = \left( \begin{array}{cccc}
a_1 & b_1 & c_1 & d_1 \\
a_2 & b_2 & c_2 & d_2 \\
. & . & . & . \\
. & . & . & . \\
a_{\bf M} & b_{\bf M} & c_{\bf M} & d_{\bf M} \end{array} \right)\]
where each row contains the coefficients of a spline interpolating
two points in a single trajectory. Then the {\bf M}x({\bf V}+1) 
matrix product {\bf E} = {\bf C}*{\bf T} contains in each row the 
values of a cubic between two data points, for all the {\bf M} 
trajectories ({\em eulerian view}). In a similar way on can consider
a {\em lagrangian view} for computing the values of all the cubics
in a single trajectory. It can be easily shown that the total number
of operations for computing all the values along each trajectory is
of order {\bf N}x{\bf M}x({\bf V}+1).

\section{Computing values of splines}
For the computation of the values we have used the cluster of our
company with multisessions of MATLAB package as parallel environment.
It is fundamental for this step the improvement of performances due
to the usage of LAPACK level 3 Blas routines incorporated in Matlab [5].
Another feature of this method is the fact that the matrix {\bf T} is
constant, hence it is computed only once, requires a small memory 
allocation so its values can be stored permanently in the cache.
With {\bf P}, number of processes, divisor of 3{\bf N}, total number of
two-points groups, the method used has been the distribution of 
3{\bf N}/{\bf P} matrices {\bf C} to every process.

The performances of multiprocess products show a quite linear speedup
respect the {\bf P} variable and a total computation time of order 
{\bf N}x{\bf M}; increasing the value of {\bf M} or {\bf N} for a better 
resolution, the time spent on computation doesn't change if the value of
processes is increased ({\em Gustafson Law} [6]) (see Fig. 5).

\begin{figure}[]
\begin{center}
\includegraphics[height=6.0cm,width=6.0cm]{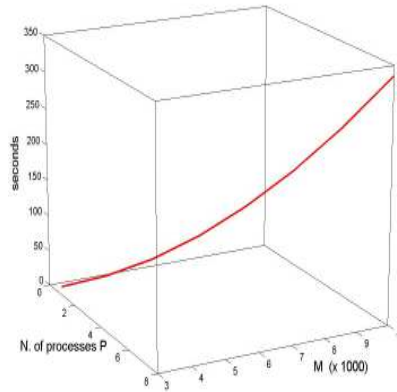}
\end{center}
\caption{Total execution time for matrices products in the case 3{\bf N}=
4200{\bf P}.}\label{fig:Perf1}
\end{figure}

\section{Conclusions}
These techniques have supplied good results for improving performances 
of post-processing phase in CFD simulations. Further work is planned for
implementing a {\em global matrix} product for the splines evaluation, with the
purpose of using the sparse matrices benefits to reduce total execution time
and memory allocation.

\end{document}